\pdfoutput=1
\RequirePackage{ifpdf}
\ifpdf 
\documentclass[pdftex]{sigma}
\else
\documentclass{sigma}
\fi

\newtheorem{Theorem}{Theorem}[section]
\theoremstyle{definition}

\numberwithin{equation}{section}

\begin{document}
\allowdisplaybreaks

\newcommand{\arXivNumber}{2010.10321}

\renewcommand{\PaperNumber}{140}

\FirstPageHeading

\ShortArticleName{An Explicit Example of Polynomials Orthogonal on the Unit Circle}

\ArticleName{An Explicit Example of Polynomials Orthogonal\\ on the Unit Circle with a Dense Point Spectrum\\ Generated by a Geometric Distribution}

\Author{Alexei ZHEDANOV}

\AuthorNameForHeading{A.~Zhedanov}

\Address{School of Mathematics, Renmin University of China, Beijing 100872, China}
\Email{\href{mailto:zhedanov@ruc.edu.cn}{zhedanov@ruc.edu.cn}}

\ArticleDates{Received November 02, 2020, in final form December 19, 2020; Published online December 21, 2020}

\Abstract{We present a new explicit family of polynomials orthogonal on the unit circle with a dense point spectrum. This family is expressed in terms of $q$-hypergeometric function of type ${_2}\phi_1$. The orthogonality measure is the wrapped geometric distribution. Some ``classical'' properties of the above polynomials are presented.}


\Classification{33D45; 42C05}

\section{Introduction}

Let $\Phi_n(z)$ be monic polynomials $\Phi_n(z)=z^n + O\big(z^{n-1}\big)$ defined through the recurrence relation~\cite{Simon}
\begin{gather*}
\Phi_{n+1}(z) =z \Phi_n(z) - \bar a_n \Phi^*_n(z), \qquad \Phi_0(z) =1, 
\end{gather*}
where
\[
\Phi_n^*(z) =z^n \bar \Phi_n(1/z)
\]
and where $\bar \Phi_n(z)$ means complex conjugation of expansion coefficients of the polynomial $\Phi_n(z)$. The recursion parameters
\[
a_n = -\bar \Phi_{n+1}(0)
\]
are called the Verblunsky (sometimes also reflection, Schur etc.) parameters~\cite{Simon}.

Under the condition
\begin{gather}
|a_n|<1, \qquad n=0,1,2,\dots \label{a_cond} \end{gather}
the polynomials $\Phi_n(z)$ are orthogonal on the unit circle with respect to a positive measure ${\rm d}\sigma(\theta)$
\begin{gather}
\int_0^{2 \pi} \Phi_n\big({\rm e}^{{\rm i} \theta}\big) \bar \Phi_m\big({\rm e}^{-{\rm i}\theta}\big){\rm d} \sigma(\theta) = h_n \delta_{nm}, \label{ort_Phi} \end{gather}
where
\begin{gather}
h_n = \big(1-|a_0|^2\big)\big(1-|a_1|^2\big) \cdots \big(1-|a_{n-1}|^2\big) \label{h_n} \end{gather}
are normalization constants (which are nonzero due to condition~\eqref{a_cond}). In this case $\Phi_n(z)$ are called the orthogonal polynomials on the unit circle (OPUC).

Note that orthogonality relation~\eqref{ort_Phi} is equivalent to conditions~\cite{Simon}
\begin{gather}
I_{nj} \equiv \int_0^{2 \pi} \Phi_n\big({\rm e}^{{\rm i} \theta}\big) {\rm e}^{-{\rm i}j\theta}\, {\rm d} \sigma(\theta) = h_n \delta_{nj}, \qquad j=0,1,2,\dots, n. \label{ort_j} \end{gather}

Equivalently, OPUC $\Phi_n(z)$ can be constructed in terms of trigonometric moments $\sigma_n$. The latter are defined as
\begin{gather*}
\sigma_n = \int_0^{2 \pi} {\rm e}^{{\rm i} n \theta} \,{\rm d} \sigma(\theta), \qquad n=0, \pm 1, \pm 2, \dots. 
\end{gather*}
Then polynomials $\Phi_n(z)$ have the explicit expression
\begin{gather*}
 \Phi_n(z)=(\Delta_n)^{-1} \left |
\begin{matrix} \sigma_0 & \sigma_1 & \dots & \sigma_n \\ \sigma_{-1} & \sigma_0 &
\dots & \sigma_{n-1} \\ \dots & \dots & \dots & \dots\\
\sigma_{1-n}& \sigma_{2-n}& \dots & \sigma_1\\ 1& z & \dots & z^n
\end{matrix} \right |, 
\end{gather*}
where
\begin{gather*} \Delta_n= \left |
\begin{matrix} \sigma_0 & \sigma_{1} & \dots &
\sigma_{n-1}\\ \sigma_{-1}& \sigma_0 & \dots & \sigma_{n-2}\\ \dots & \dots & \dots & \dots\\
\sigma_{1-n} & \sigma_{2-n} & \dots & \sigma_0 \end{matrix} \right| 
\end{gather*}
are Toeplitz determinants which are all positive $\Delta_n>0$, $n=0,1,2,\dots$.
Note the symmetry property of the trigonometric moments
\begin{gather}\label{sym_sigma}
\sigma_{-n} = \bar \sigma_n.
\end{gather}
Explicit examples of polynomials orthogonal on unit circle are very interesting from different point view. By ``explicit examples'' we mean that all main objects: the parameters~$a_n$, the moments~$\sigma_n$, the measure $\sigma(\theta)$ and the polynomials themselves $\Phi_n(z)$ have explicit expressions in terms of special functions. Usually, in most known explicit examples the parameters~$a_n$ are given by elementary functions of~$n$ while the OPUC~$\Phi_n(z)$ are expressed in terms of hypergeometric functions (either ordinary or basic). A list of known explicit examples can be found, e.g., in Simon's monograph~\cite{Simon}.

In \cite{Tsu_Zhe,Zhe_cndn} new explicit examples of OPUC were presented. In these examples polyno\-mials~$\Phi_n(z)$ are expressed in terms of elliptic hypergeometric function ${_3}E_2(z)$ while the moments~$\sigma_n$ and the recurrence parameters $a_n$ have simple expressions in terms of elliptic functions. The most interesting property of the OPUC of these examples is that they are orthogonal on the unit circle with respect to a dense point measure. This means that the function~$\sigma(\theta)$ is a step function with infinitely many points~$\theta_s$ of jumps, and these points are dense on the interval~$[0, 2 \pi]$. In terms of the distribution function this can be presented as
\begin{gather*}
\rho(\theta) = \sum_{s=-\infty}^{\infty} M_s \delta(\theta-\theta_s), 
\end{gather*}
where $\rho(\theta)$ is a distribution defined as ${\rm d} \sigma(\theta) = \rho(\theta) {\rm d} \theta$, $\delta(\theta)$ is the Dirac delta function and~$M_s$ are concentrated masses located at points of jumps~$\theta_s$. The spectral points $z_s=\exp({\rm i} \theta_s)$ are dense on the unit circle.

Then orthogonality relation \eqref{ort_Phi} can then be presented as
\[
\sum_{s=-\infty}^{\infty} M_s \Phi_n\big({\rm e}^{{\rm i} \theta_s}\big) \bar \Phi_m\big({\rm e}^{-{\rm i} \theta_s}\big) = h_n \delta_{nm}. \]

From general considerations (see, e.g.,~\cite{Simon}) it follows that polynomials orthogonal with respect to such dense point measures are rather generic if one assumes some natural restrictions upon behavior of the recurrence parameters~$a_n$. On the other hand, such measures are very important from physical point of view, because they correspond to the phenomenon of the Anderson localization \cite{50,Simon}.

Usually examples of OPUC with dense point spectrum are related to sequences of the para\-me\-ters~$a_n$ which behave (quasi) stochastically inside the interval $|a_n|<1$ \cite{Simon}. OPUC in~\cite{Tsu_Zhe,Zhe_cndn}  provide perhaps the first known examples of a pure point dense measure on the unit circle where both the coefficients~$a_n$ and the moments~$\sigma_n$ are given explicitly by analytic functions in~$n$.

In this paper we propose a much simpler explicit example of polynomials orthogonal on the unit with respect to a (wrapped) geometric distribution which is dense on the unit circle. Polynomials themselves are expressed in terms of basic hypergeometric function ${_2}\phi_1(q;z)$ with $|q|=1$.

\section{Wrapped geometric distribution and corresponding OPUC}

Let $q$ be a fixed point belonging to the unit circle $|q|=1$ and not a root of unity, i.e., we demand that $q^n \ne 1$ for all natural integer $n=0,1, \dots$. Choose a real parameter $p$ within the unit interval $0<p<1$. Define the measure on the unit circle as
\begin{gather}
\rho(\theta) = (1-p)\sum_{s=0}^{\infty} p^s \delta(\theta - s \chi), \label{geom_circle} \end{gather}
where $\chi$ is a fixed {\it irrational} parameter $0< \chi< 1$ such that
\begin{gather}
q=\exp(2 \pi {\rm i} \chi). \label{q_chi} \end{gather}
Irrationality of $\chi$ means that the set of points $z_s=q^s$, $s=0,1,\dots$ (i.e., the location of jumps of the measure) is dense on the unit circle. The weights (i.e., the concentrated masses)~$w_s$ at the points~$z_s$ form the geometric sequence: $w_s=p^s$, $s=0,1,2,\dots$.

Corresponding trigonometric moments are
\begin{gather}
\sigma_n = \sum_{s=0}^{\infty} z_s^n w_s = (1-p) \sum_{s=0}^{\infty} q^{sn} p^s = \frac{1-p}{1-pq^n}, \qquad n=0, \pm 1, \pm 2, \dots . \label{mom_geom} \end{gather}
Note that the measure~\eqref{geom_circle} can be interpreted as the wrapped geometric distribution on the unit circle (see, e.g., \cite{Jacob, MJ} for definition and discussion of wrapped distributions on the unit circle).

Relation $\sigma_0=1$ means that the measure~\eqref{geom_circle} has the standard normalization condition.

One can present expression \eqref{mom_geom} as
\begin{gather}
\sigma_n = \frac{(p;q)_n}{(pq;q)_n}, \label{sigma_poch} \end{gather}
where the $q$-shifted factorial ($q$-Pochhammer symbol) is defined as~\cite{GR,KS} $(a;q)_0=1$ and
\begin{gather*}
(a;q)_n = (1-a)(1-aq) \cdots \big(1-aq^{n-1}\big) 
\end{gather*}
for positive $n=1,2,\dots$ and
\begin{gather*}
(a;q)_n = \frac{1}{\big(aq^{n};q\big)_{-n}} 
\end{gather*}
for negative $n=-1,-2,\dots$.

It is known that the Laurent biorthogonal Pastro polynomials $P(z;a,b)$~\cite{Pas} depending on two arbitrary parameters $a$, $b$ can be uniquely defined through their moments~\cite{VZ}
\begin{gather}
\sigma_n = \frac{(a;q)_n}{(b;q)_n}, \qquad n=0, \pm 1, \pm 2, \dots. \label{sigma_Pastro}
\end{gather}
Explicitly, these polynomials are given by \cite{Pas, VZ}
\begin{gather}
P(z;a,b) = \mu_n \, {_2}\phi_1 \left( { q^{-n}, b \atop a q^{1-n}} ; qz \right), \label{Pastro_P} \end{gather}
where $\mu_n$ is an appropriate normalization factor to fulfill the condition $P_n(z) = z^n + O\big(z^{n-1}\big)$. The definition and notation of the basic hypergeometric function ${_m}\phi_n(z)$ is standard (see, e.g., \cite{GR,KS}). For example, in the special case $m=n+1$ we have the expressions
\[
{_{m+1}}\phi_m \left( { a_1, a_2 \dots, a_{m+1}\atop b_1, b_2 \dots, b_m }; z \right) = \sum_{s=0}^{\infty} \frac{(a_1;q)_s (a_2;q)_s \cdots (a_{m+1};q)_s}{(q;q)_s(b_1;q)_s (b_2;q)_s \cdots (b_{m};q)_s}z^s.
\]

Note the Laurent biorthogonal polynomials (LBP) can be considered as a generalization of the OPUC. Their main distinction from OPUC is that the moments~$\sigma_n$ do not satisfy, in general, the symmetry condition~\eqref{sym_sigma}. The LBP can also be characterized by the three-term recurrence relation of~$R_I$ type~\cite{HR, Zh_LBP}
\begin{gather}
P_{n+1}(z) + g_n P_n(z) =z  ( P_n(z) + d_n P_{n-1}(z) ) , \qquad P_0=1, \qquad P_{-1}=0 \label{3term_P}
\end{gather}
with some recurrence coefficients $g_n$, $d_n$.

In contrast to the case of orthogonal polynomials, the recurrence relation \eqref{3term_P} can be presented in the form of the generalized eigenvalue problem \cite{Zh_GEVP}
\begin{gather}
J_1 {\mathbf P}(z) = z J_2 {\mathbf P}(z), \label{GEVP_L} \end{gather}
where $J_1$, $J_2$ are upper and lower bidiagonal matrices acting on the vector
\[\mathbf P(z) = (P_0(z), P_1(z), \dots).\]

Comparing expressions \eqref{sigma_Pastro} and \eqref{sigma_poch} we can conclude that the OPUC corresponding to the wrappedd geometric distribution are special case of the Pastro polynomials with $a=p$, $b=qp$.

This allows us to present the main result

\begin{Theorem}
The polynomials $\Phi_n(z)$ orthogonal on the unit circle with respect to the wrapped geometric distribution \eqref{geom_circle} have the explicit expression
\begin{gather} \Phi_n(z) =
\mu_n \: {_2}{\phi}_1\left( {q^{-n}, pq \atop pq^{1-n}};zq \right), \label{Pastro_geom}
\end{gather} where
\begin{gather}
\mu_n = q^{-n}  \frac{(q;q)_n\big(pq^{1-n};q\big)_n}{\big(q^{-n};q\big)_n(pq;q)_n}= p^n \frac{\big(p^{-1};q\big)_n}{(qp;q)_n}.\label{mu_n}
\end{gather}
\end{Theorem}

One can directly check that the polynomials \eqref{Pastro_geom} satisfy orthogonality relations \eqref{ort_j}. Indeed, one has
\begin{gather*}
I_{nj}= (1-p) \sum_{s=0}^{\infty} \Phi_n\big(q^s\big) q^{-sj} p^s = (1-p) \mu_n \sum_{s=0}^{\infty} \sum_{k=0}^n \frac{\big(q^{-n};q\big)_k (pq;q)_k}{(q;q)_k \big(pq^{1-n};q\big)_k} q^{(1+s)k} q^{-js} p^s . 
\end{gather*}
Performing summation over $s$ we get
\begin{gather*}
 I_{nj}= \frac{(1-p) \mu_n}{1-pq^{-j}}   \sum_{k=0}^n \frac{\big(q^{-n};q\big)_k (pq;q)_k \big(pq^{-j};q\big)_k }{(q;q)_k \big(pq^{1-n}\big)_k \big(pq^{1-j};q\big)_k} q^k =
 \frac{(1-p) \mu_n}{1-pq^{-j}} \, {_3}{\phi}_2 \left( {q^{-n}, pq, pq^{-j} \atop pq^{1-n},pq^{1-j}};q \right).
\end{gather*}
The above expression can be simplified by the $q$-Saalsch\"utz formula \cite{GR,KS}:
\begin{gather*}
{_3}{\phi}_2 \left( {q^{-n}, a, b \atop c, abc^{-1} q^{1-n}};q \right) = \frac{(c/a;q)_n (c/b;q)_n}{(c;q)_n (c/(ab);q)_n}.
\end{gather*}
We thus have
\begin{gather}
I_{nj} = \frac{(1-p) \mu_n}{1-pq^{-j}}   \frac{\big(q^{-n};q\big)_n\big(q^{j+1-n};q\big)_n}{\big(pq^{1-n};q\big)_n \big(p^{-1} q^{j-n};q\big)_n}.\label{I_3}
\end{gather}
The factor $\big(q^{j+1-n};q\big)_n$ in~\eqref{I_3} becomes zero when $j=0,1,\dots, n-1$ and hence
\[
I_{nj}=0, \qquad j=0,1,\dots, n-1,
\]
which is equivalent to orthogonality relation \eqref{ort_j}. It remains to show that $I_{nn}>0$. After simple calculations one can arrive at the expression
\begin{gather}
I_{nn}=h_n=\frac{|(q;q)_n|^2}{|(pq;q)_n|^2}   p^n
\label{h} \end{gather}
from which it is clear that $h_n>0$ for all $n=0,1,2,\dots$ due to condition $q^n \ne 1$.

Explicit expression for the the recurrence parameters $a_n$ follows from \eqref{Pastro_geom} and \eqref{mu_n}:
\begin{gather*}
{\bar a}_{n-1} = - \Phi_n(0) = -\mu_n = - p^n   \frac{\big(p^{-1};q\big)_n}{(qp;q)_n}. 
\end{gather*}
For the square of absolute values we have rather simple expression
\begin{gather}
|a_{n-1}|^2 = {\bar a}_{n-1} a_{n-1} = \frac{(1-p)^2}{1+p^2 -p \big(q^n+q^{-n}\big)} = \frac{1}{1+\beta \sin^2 (\chi \pi n)}, \label{abs_a}
\end{gather}
where
\begin{gather*}
\beta= \frac{4 p}{(1-p)^2} 
\end{gather*}
and where the parameter $\chi$ is the same as in~\eqref{q_chi}.

It is seen from \eqref{abs_a} that the values $|a_n|$ oscillate inside the interval
\begin{gather*}
\frac{1-p}{1+p} < |a_n| <1, \qquad n=0,1,\dots . 
\end{gather*}
Because of irrationality of $\chi$ the absolute value $|a_n|$ achieves the boundaries of this interval with any prescribed accuracy (never achieving exact boundary values). Note that $a_{-1}=-1$ which corresponds to the standard initial conditions for OPUC~\cite{Ger,Simon}. It is easily verified that expression~\eqref{h} for~$h_n$ agrees with relation~\eqref{h_n}.

The OPUC \eqref{Pastro_geom} can be considered as $|q|=1$ analogs of the OPUC introduced by Askey in~\cite{Ask_Sz} (see also~\cite{Costa} for more general OPUC of Askey's type).

\section[Classicality of the polynomials Phin(z)]{``Classicality'' of the polynomials $\boldsymbol{\Phi_n(z)}$}

The OPUC \eqref{Pastro_geom} possess ``classical'' properties which make them similar to classical orthogonal polynomials.

First of all, they satisfy the three-term recurrence relation \eqref{3term_P} where the recurrence coefficients are
\begin{gather*}
g_n = \frac{q^n-p}{1-p q^{n+1}}, \qquad d_n = - \frac{p(1-q^n)^2}{(1-p q^{n})(1-p q^{n+1})}. 
\end{gather*}
Moreover, the polynomial $\Phi_n(z)$ possess a remarkable {\it duality} property. Indeed, one can rewrite polynomials $\Phi_n(z)$
in a different form
\begin{gather} \Phi_n(z) =
p^n \frac{(q;q)_n}{(pq;q)_n} z^n \,
{_3}{\phi}_2\left( {q^{-n}, p^{-1}, z^{-1} \atop q,0}; q
\right), \label{32_form}
\end{gather} which can be obtained from
\eqref{Pastro_geom} by standard transformation formulas \cite{GR,KS}.

From this formula the duality property
\begin{gather}
A_s \Phi_s(q^n) = A_n  \Phi_n(q^s) \label{dual}
\end{gather}
follows, where
\begin{gather*}
A_n = \frac{(pq;q)_n}{(q;q)_n}   p^{-n}. 
\end{gather*}
This property resembles corresponding duality properties for the classical orthogonal polynomials from the Askey scheme \cite{BI,Leonard,Ter}. The main difference is that the polynomials $\Phi_n(z)$ satisfy the generalized eigenvalue problem~\eqref{GEVP_L} instead of the ordinary eigenvalue problem for orthogonal polynomials.

From the duality property one can derive the second-order $q$-difference equation
\begin{gather}
B_{s+1} \Phi_n\big(q^{s+1}\big) + g_s \Phi_n\big(q^{s}\big) = q^n \big(\Phi_n\big(q^{s}\big) + B_s^{-1} d_s \Phi_n\big(q^{s-1}\big) \big), \label{dfr_eq}
\end{gather}
where
\[
B_s = \frac{A_{s-1}}{A_s} = \frac{p\big(1-q^s\big)}{1-pq^s}.
\]
Equation \eqref{dfr_eq} can also be presented in the form of the generalized eigenvalue problem
\begin{gather}
L \Phi_n(z) = q^n M \Phi_n(z),
\label{LM_GEVP} \end{gather}
where the first-order $q$-difference operators $L$, $M$ act on the argument $z$ of the polynomials.

Relations \eqref{GEVP_L} and \eqref{LM_GEVP} mean that the polynomials possess the {\it bispectrality} property: they satisfy simultaneously two GEVP. Concerning definition and general theory of bispectrality see, e.g.,~\cite{Grunbaum}. For orthogonal polynomials from the Askey scheme this property is well known~\cite{KS}. For biorthogonal polynomials and rational functions the bispectrality is known for some special fa\-milies. The most general from them are elliptic biorthogonal functions \cite{SZ}. However the general theory of bispectrality for systems satisfying GEVP is not yet developed (see, e.g., \cite{VZ_Hahn1,VZ_Hahn2} for algebraic description of bispectrality on the ``lowest'' level of hypergeometric functions~${_3}F_2(1)$).

The duality property implies that for $z=q^s$, $s=0,1,2,\dots$
the hypergeometric function in~\eqref{32_form} reduces to a~polynomial of degree~$s$ of the argument~$q^{-n}$.

It is well known (see, e.g.,~\cite{Simon}) that if $z_0$ is a point
on the unit circle corresponding to a~concentrated mass $M_0$ then
the relation
\begin{gather*} \sum_{n=0}^{\infty} \frac{|\Phi_n(z_0)|^2}{h_n} = 1/M_0, \label{bound_sum}
\end{gather*}
holds, where the normalization coefficient $h_n$ is defined in~\eqref{h_n}.

In our case this means that for every spectral point $z_s=q^s$, $s=0,1,\dots$ there exists the identity
\begin{gather}
\sum_{n=0}^{\infty} \frac{\big|\Phi_n\big(q^s\big)\big|^2}{h_n} = {M_s}^{-1} = \frac{p^{-s}}{1-p}. \label{sum_discr}
\end{gather}
Identity \eqref{sum_discr} follows easily from the duality property~\eqref{dual} and from orthogonality relation.

So far, we have considered the case when $q$ is not a root of unity. If, otherwise, $q$ is a primitive root of unity
\begin{gather*}
q=\exp\left( \frac{2 \pi {\rm i} M}{N} \right) 
\end{gather*}
with coprime integers $M$, $N$, then there are only $N$ distinct mass points on
the unit circle located at $z_s=q^s$, $ s=0,1,2,\dots,
N-1$. In this case the polynomials $\Phi_n(z)$ are orthogonal on vertices of a regular
$N$-gon with respect to the finite wrapped geometric distribution:
\begin{gather*} \sum_{s=0}^{N-1} \Phi_n\big(q^s\big) \bar \Phi_m\big(q^{-s}\big) \big(1-p^N\big)
p^s = h_n  \delta_{nm}, \qquad n,m=0,1,\dots, N-1 .
\end{gather*}
See \cite{Zhe2} for other explicit examples of
polynomials orthogonal on the vertices of regular polygons.

\section{Concluding remarks}

In contrast to examples of OPUC obtained in \cite{Zhe_cndn}, the polynomials~\eqref{Pastro_P} have non-real moments~$\sigma_n$ and hence the coefficients~$a_n$ are non-real as well. This means that it is impossible to associate with OPUC~\eqref{Pastro_P} polynomials orthogonal on an interval of the real line. In~\cite{Zhe_cndn} explicit examples of polynomials orthogonal with dense point spectrum on an interval were presented using standard Szeg\H{o} mapping from OPUC to an interval of the real line. We mention also examples of OPUC and ordinary orthogonal polynomials with dense point spectrum presented in \cite{Magnus_snul, Magnus_semi}.

The OPUC \eqref{Pastro_P} allow a trivial modification which shifts all spectral points on the unit circle on the same constant angle $\varphi$, i.e., we can consider the same weights $w_s=p^s (1-p)$ located at the points
\[
\theta_s = 2 \pi \chi s + \varphi, \qquad s=0,1,2,\dots .
\]
Equivalently, this means that the new spectral points will be $\tilde z_s = {\rm e}^{{\rm i} \varphi} q^s$, $s=0,1,2,\dots$.

Such transformation is equivalent to a simple rotation of the argument of OPUC \cite{Ger,Simon}:
\[
\tilde \Phi_n(z) = {\rm e}^{-{\rm i} \varphi n} \Phi_n\big({\rm e}^{{\rm i} \varphi } z\big).
\]
Another modification of the OPUC \eqref{Pastro_P} is more substantional. It leads to Laurent biorthogonal polynomials orthogonal on the unit circle with dense point measure.

Indeed, assume that the spectral points on the unit circle are the same: $z_s=q^s, \, s=0,1,2, \dots$. Take the weights:
\begin{gather}
w_s = p^s \frac{\big(q^k;q\big)_s}{(q;q)_s}, \qquad 0<p<1, \qquad k=1,2,3,\dots.
\label{w_s_k} \end{gather}
For $k=1$ we return to the case of the wrapped geometric distribution. For $k>1$ the moments are
\begin{gather}
\sigma_n = \sum_{s=0}^{\infty} \frac{\big(q^k;q\big)_s}{(q;q)_s} p^s q^{sn}. \label{sigma_k} \end{gather}
By $q$-binomial theorem \cite{GR,KS} the above sum is simplified to
\begin{gather*}
\sigma_n = \frac{\big(q^k pq^n;q\big)_{\infty}}{\big(pq^n;q\big)_{\infty}} = \frac{1}{\big(pq^n;q\big)_k} = \frac{(p;q)_n}{(p;q)_k \big(pq^k;q\big)_n}. \label{sig_k}
\end{gather*}

\begin{remark*}
Usually, the convergence problem for $q$-series like \eqref{sigma_k} with $|q|=1$ is highly nontrivial (see, e.g.,~\cite{LubSaff}). In our case however this problem does not appear because for integer~$k$ there is cancellation of almost all terms (apart of a finite number of initial ones) in denominators of the coefficients in~\eqref{sigma_k}. Hence the convergence for $0<p<1$ still takes place.
\end{remark*}

For fixed $k$ the moments $\sigma_n$ coincide (up to a constant factor) with the moments \eqref{sigma_Pastro} for the Pastro polynomials with $a=p$, $b=pq^k$. Hence from \eqref{Pastro_P} we have explicit expression for them
\begin{gather*}
P(z) = \mu_n   {_2}\phi_1 \left( { q^{-n}, pq^k \atop p q^{1-n}} ; qz \right), 
\end{gather*}
where
\[
\mu_n = \frac{p^n \big(p^{-1};q\big)_n}{\big(pq^k;q\big)_n}.
\]
These polynomials are NOT OPUC (apart from the already considered case $k=1$) because the weights~\eqref{w_s_k} are not positive and hence the moments $\sigma_n$ do not satisfy symmetry property~\eqref{sym_sigma}.

Existence of other explicit examples of OPUC with dense point spectrum is an interesting open problem.

\subsection*{Acknowledgements}

The author is indebted for F.A.~Gr\"unbaum, A.~Magnus, V.~Spiridonov, S.~Tsujimoto and L.~Vinet for discussions and for anonymous referees for valuable remarks. The author is gratefully holding Simons CRM Professorship and is funded by the National Foundation of China (Grant No.~11771015).

\pdfbookmark[1]{References}{ref}
\LastPageEnding

\end{document}